\documentclass[global]{svjour}
\usepackage{amssymb}
\usepackage{amsfonts}
\usepackage{amsmath}

\input{tcilatex}

\begin{document}

\title{Generalized Vector Quasi-Equilibrium Problems Involving Relaxed
Continuity of the Set-Valued Maps}
\author{Monica Patriche}
\institute{University of Bucharest 
\email{monica.patriche@yahoo.com}%
}
\mail{\\
University of Bucharest, Faculty of Mathematics and Computer Science, 14
Academiei Street, 010014 Bucharest, Romania}
\maketitle

\textit{University of Bucharest,} \textit{Faculty of Mathematics and
Computer Science, }

\textit{Academiei 14 Street,010014Bucharest,Romania.}

\textit{E-mail:monica.patriche@yahoo.com }\bigskip

\textbf{Abstract.}{\small \ In this paper, we study }t{\small he existence
of solutions for generalized vector quasi-equilibrium problems. Firstly, we
prove that in the case of Banach spaces, the assumptions of continuity over
correspondences can be weakened. The theoretical analysis is based on
fixed-point theorems. Secondly, we establish new equilibrium theorems as
applications of the KKM principle.}


\section{Introduction}

The vector equilibrium problem is a unified model of several other problems,
for instance, vector variational inequalities, vector optimization problems
or Debreu-type equilibrium problems. This explains the increasing attention
paid to it by many researchers nowadays. For further information on this
topic, the reader is referred to the following list of selected
publications: [1-5],\cite{far},\cite{far2},\cite{hou}, \cite{lin},[14-17].

We report new results concerning the existence of solutions for the
generalized vector quasi-equilibrium problems with correspondences
fulfilling weak-continuity assumptions. Firstly, {\normalsize our
theoretical analysis is based on fixed-point theorems. In the case of Banach
spaces, the assumption of continuity over correspondences can be weakened to
almost lower semicontinuity. The methodology of the proofs relies on an
approximation technique. We also establish a new equilibrium result in the
case of lower semicontinuity of the correspondences. }We work on regular
vector topological spaces (dropping the local convexity). We also drop the
property of properly $C-$ quasi convexity of $f(x,y,\cdot ).${\normalsize \
Our study is characterized by the refinement of the hypotheses concerning
the equilibrium problems. We underline its novelty and the improvement of
the existent results obtained so far in literature. Secondly,} we obtained
new equilibrium theorems by applying the KKM principle.

The rest of the paper is organized as follows. Section 2 contains
preliminaries and notations. Sections 3 reveals the existence of solutions
for a generalized vector quasi-equilibrium problem. Several applications of
the KKM principle to the equilibrium problems are established in Section 4.
Section 5 presents the conclusions of our research.

\section{Preliminaries and notation}

For the reader's convenience, we present several properties of the
correspondences which are used in our proofs.

Let $X$, $Y$ be topological spaces and $T:X\rightrightarrows Y$ be a
correspondence. $T$ is said to be \textit{upper semicontinuous} if for each $%
x\in X$ and each open set $V$ in $Y$ with $T(x)\subset V$, there exists an
open neighborhood $U$ of $x$ in $X$ such that $T(x)\subset V$ for each $y\in
U$. $T$ is said to be \textit{lower semicontinuous} if for each x$\in X$ and
each open set $V$ in $Y$ with $T(x)\cap V\neq \emptyset $, there exists an
open neighborhood $U$ of $x$ in $X$ such that $T(y)\cap V\neq \emptyset $
for each $y\in U$. Equivalently, $T$ is \textit{lower semicontinuous} if for
each closed set $V$ in $Y,$ $\{x\in X:T(x)\subset V\}$ is closed. For each $%
x\in X,$ the set $T(x)$ is called the \textit{upper section} of $T$ at $x.$
For each $y\in Y,$ the set $T^{-1}(y):=\{x\in X:y\in T(x)\}$ is called the 
\textit{lower section} of $T$ at $y$.

The graph of $T:X\rightrightarrows Y$ is the set Gr$(T)=\{(x,y)\in X\times
Y:y\in T(x)\}.$

The correspondence $T:X\rightrightarrows Y$ is said to have \textit{a closed
graph} if the set Gr$(T)$ is closed in $X\times Y.$ If $T$ has a closed
graph, then, the lower and the upper sections are closed.

The correspondence $\overline{T}$ is defined by $\overline{T}(x)=\{y\in
Y:(x,y)\in $cl$_{X\times Y}$Gr$(T)\}$ (the set cl$_{X\times Y}$Gr$(T)$ is
called \textit{the adherence of the graph of }$T$)$.$ It is easy to see that
cl$T(x)\subset \overline{T}(x)$ for each $x\in X.$

Lemma 1 concerns the continuity of correspondences. It will be also crucial
in our proofs.

\begin{lemma}
(see \cite{yuan}). \textit{Let }$X$\textit{\ and }$Y$\textit{\ be two
topological spaces and let }$A$\textit{\ be a closed subset of }$X.$\textit{%
\ Suppose }$T_{1}:X\rightrightarrows Y$\textit{, }$T_{2}:X\rightrightarrows
Y $\textit{\ are lower semicontinuous correspondences such that }$%
T_{2}(x)\subset T_{1}(x)$\textit{\ for all }$x\in A.$\textit{\ Then the
correspondence }$T:X\rightrightarrows Y$\textit{\ defined by }$T\mathit{(z)=}%
\left\{ 
\begin{array}{c}
T_{1}(x)\text{, \ \ \ \ \ \ \ if }x\notin A\text{, } \\ 
T_{2}(x)\text{, \ \ \ \ \ \ \ \ if }x\in A%
\end{array}%
\right. $\textit{is also lower semicontinuous.}
\end{lemma}

Let $X$ be a nonempty convex subset of a topological vector space\textit{\ }$%
E,$ $Z$ be a real topological vector space, $Y$ be a subset of $Z$ and $C$
be a pointed closed convex cone in $Z$ with its interior int$C\neq \emptyset
.$ Let $T:X\rightrightarrows Z$ be a correspondence with nonempty values. $T$
is said to be \textit{properly }$C-$\textit{quasi-convex on} $X$, if for any 
$x_{1},x_{2}\in X$ and $\lambda \in \lbrack 0,1],$ either $T(x_{1})\subset
T(\lambda x_{1}+(1-\lambda )x_{2})+C$ or $T(x_{2})\subset T(\lambda
x_{1}+(1-\lambda )x_{2})+C.$

\section{Main results}

We start with the presentation of the problem we approach.\medskip

Let $X$ be a nonempty subset of a topological vector space $E$ and let $Y,Z$
be topological vector spaces. Let $C:X\rightrightarrows Z$ be a
correspondence such that, for each $x\in X,$ $C(x)$ is a pointed closed
convex cone with int$C(x)\neq \emptyset .$ Let $K:X\rightrightarrows X,$ $%
T:X\rightrightarrows Y$and $f:X\times Y\times X\rightrightarrows Z$ be
correspondences with nonempty values. We consider the following generalized
vector quasi-equilibrium problem:

(GVQEP 1): Find $(x^{\ast },y^{\ast })\in X\times Y$ such that $x^{\ast }\in 
\overline{K}(x^{\ast }),$ $y^{\ast }\in \overline{T}(x^{\ast })$ and $%
f(x^{\ast },y^{\ast },u)\nsubseteq -$int$C(x^{\ast })$ for each $u\in
K(x^{\ast })$.$\medskip $

Firstly, we study a generalized vector quasi-equilibrium problem on Banach
spaces with the general motivation of relaxing the continuity assumptions of
the involved correspondences. The almost lower semicontinuous correspondences%
\textit{\ }defined by Deutsch and Kenderov in \cite{deu} enjoy the very
interesting property of having $\varepsilon -$approximate selections. This
leads to the fact that finding solutions for GVQEP by using this property is
becoming increasing crucial.\smallskip

Now, we are going to present the almost lower semicontinuous
correspondences.\smallskip

Let $X$\ be a topological space and $Y$\ be a normed linear space.\textit{\ }%
The correspondence $T:X\rightrightarrows Y$ is said to be \textit{almost
lower semicontinuous (a.l.s.c.)} \textit{at} $x\in X$ (see \cite{deu}), if,
for any $\varepsilon >0,$ there exists a neighborhood $U(x)$ of $x$ such
that $\tbigcap\limits_{z\in U(x)}B(T(z);\varepsilon )\neq \emptyset .$

T is \textit{almost lower semicontinuous} if it is a.l.s.c. at each $x\in X$%
.\smallskip

Deutsch and Kenderov \cite{deu} established the following characterization
of a.l.s.c. correspondences.

\begin{lemma}
\textit{(Deutsch and Kenderov, \cite{deu}) Let }$X$\textit{\ be a
paracompact topological space, }$Y$\textit{\ be a normed vector space and }$%
T:X\rightrightarrows Y$\textit{\ be a correspondence having convex values.
Then, }$\mathit{T}$\textit{\ is a.l.s.c. if and only if, for each }$%
\varepsilon >0$\textit{, }$T$\textit{\ admits a continuous }$\varepsilon -$%
\textit{approximate selection f; that is, }$f:X\rightrightarrows Y$\textit{\
is a continuous single-valued function, such that }$f(x)\in
B(T(x);\varepsilon )$\textit{\ for each }$x\in X.\medskip $
\end{lemma}

Let $X,Y$\ and $Z$\ be topological spaces and let $K:X\rightrightarrows X,$ $%
C:X\rightrightarrows Z$\ and $f:X\times Y\times X\rightrightarrows Z$\ be
correspondences with nonempty values.\textit{\ }Let $P,G:X\times
Y\rightrightarrows X$ be defined by

$P(x,y)=\{u\in X:f(x,y,u)\subseteq -C(x)\}$ for each $(x,y)\in X\times Y$ and

$G(x,y)=K(x)\cap P(x,y)$ for each $(x,y)\in X\times Y.$

We introduce the following notion.

\begin{definition}
\textit{We say that }$f$\textit{\ has the property }$\mathcal{P}$\textit{\
with respect to }$K$\textit{\ on }$X\times Y$ \textit{if }$G$ \textit{is
almost lower semicontinuous on }$X\times Y.$
\end{definition}

We note that $f$ satisfies the following property:

$(\mathcal{P})$ \ \textit{for each }$\varepsilon >0$\textit{\ and }$(x,y)\in
X\times Y,$\textit{\ there exist }$u\in X$\textit{\ and a neighborhood }$%
U((x,y))$\textit{\ in }$X\times Y$ \textit{such that, for each }$(x^{\prime
},y^{\prime })\in U((x,y)),$\textit{\ there exists }$u^{\prime }(x^{\prime
},y^{\prime })\in K(x^{\prime })\cap B(u;\varepsilon )$\textit{\ with the
property that }$f(x^{\prime },y^{\prime },u^{\prime }(x^{\prime },y^{\prime
}))\subset -C(x^{\prime })$\textit{.}\smallskip

The next theorem establishes the existence of solutions for the considered
generalized vector quasi-equilibrium problem involving correspondences with
relaxed continuity, defined on Banach spaces.

\begin{theorem}
Let $Z$ be a Banach space and l\textit{et }$X,Y$\textit{\ be nonempty
compact convex subsets of Banach spaces }$E$\textit{, respectively }$F$%
\textit{. Assume that:}
\end{theorem}

\textit{i) }$K:X\rightrightarrows X$\textit{\ is an almost lower
semicontinuous correspondence with nonempty convex values such that }$%
\overline{K}$ \textit{has nonempty lower sections; in addition, there exists
an integer }$n_{0}\in \mathbb{N}^{\ast }$\textit{\ such that }$B(K(X);\frac{1%
}{n_{0}})\subset X;$

\textit{ii) }$T:X\rightrightarrows Y$\textit{\ is an upper semicontinuous
correspondence with nonempty convex values;}

\textit{iii) }$C:X\rightrightarrows Z$\textit{\ is a correspondence such
that, for all }$x\in X,$\textit{\ }$C(x)$\textit{\ is a pointed closed
convex cone with int}$C(x)\neq \emptyset ;$

\textit{iv) }$f:X\times Y\times X\rightrightarrows Z$\textit{\ is a
correspondence with nonempty values which satisfies the following conditions:%
}

\textit{\ \ \ \ \ \ \ \ \ \ a) \ for all }$(x,y)\in X\times Y,$\textit{\ }$%
(x,y,x)\notin $cl$\{(x,y,u)\in X\times Y\times X:f(x,y,u)\subseteq -C(x)\};$

\textit{\ \ \ \ \ \ \ \ \ \ b) for each }$(x,y)\in X\times Y,$\textit{\ the
correspondence }$f(x,y,\cdot ):X\rightrightarrows Z$\textit{\ is }$(\mathit{-%
}C(x))-$\textit{quasi-convex;}

\textit{\ \ \ \ \ \ \ \ \ \ c) }$A=\{(x,y)\in X\times Y:$\textit{\ there
exists }$u\in K(x)$\textit{\ such that}$:$\textit{\ }$f(x,y,u)\subseteq
-C(x)\}$\textit{\ is open or empty;}

\ \ \ \ \ \ \ \ \ \ \textit{d) }$f$\textit{\ has the property }$\mathcal{P}$%
\textit{\ with respect to }$K$ \textit{on} $A;$

\ \ \ \ \ \ \ \ \textit{\ e) for each }$(x_{0},y_{0})\in $\textit{fr}$A$%
\textit{, there exists an open neighborhood }$U$\textit{\ of }$(x_{0},y_{0})$%
\textit{\ such that} $\bigcap\nolimits_{(x,y)\in U\cap A}(K\cap
P)(x,y)\subset \bigcap\nolimits_{(x,y)\in U}K(x);$

\ \ \ \ \ \ \ \ \textit{f)} \textit{for each }$z\in X,$ $\{(x,y)\in X\times
Y:(x,y,z)\in $cl$_{X\times Y\times X}\{(x,y,u)\in X\times Y\times
X:f(x,y,u)\subseteq -$int$C(x)$ \textit{and} $u\in K(x)\}\}$ \textit{is
nonempty.}

\textit{Then, there exists }$(x^{\ast },y^{\ast })\in X\times Y$ \textit{%
such that }$x^{\ast }\in \overline{K}(x^{\ast }),$\textit{\ }$y^{\ast }\in 
\overline{T}(x^{\ast })$\textit{\ and }$f(x^{\ast },y^{\ast },u)\nsubseteq -$%
\textit{int}$C(x^{\ast })$\textit{\ for each }$u\in K(x^{\ast }).\medskip $

\textit{Proof. }The key idea of the proof is to use the Kakutani-Ky Fan
fixed point Theorem. To this end, we need to construct the following
correspondences.

Let $P,G:X\times Y\rightrightarrows X$ be defined by

$P(x,y)=\{u\in X:f(x,y,u)\subseteq -C(x)\}$ for each $(x,y)\in X\times Y$ and

$G(x,y)=K(x)\cap P(x,y)$ for each $(x,y)\in X\times Y.$

We are going to show that there exists a pair\textit{\ }$(x^{\ast },y^{\ast
})\in X\times Y,$ such that $x^{\ast }\in \overline{K}(x^{\ast }),$ $y^{\ast
}\in \overline{T}(x^{\ast })$ and $K(x^{\ast })\cap P(x^{\ast },y^{\ast
})=\emptyset .$

We consider two cases.

Case I.

$A=\{(x,y)\in X\times Y:$\textit{\ }$P(x,y)\cap K(x)\neq \emptyset
\}=\{(x,y)\in X\times Y:$\textit{\ }$G(x,y)\neq \emptyset \}\mathit{=}%
\emptyset $.

In this case, $G(x,y)=\emptyset $ for each $(x,y)\in X\times Y.$

Firstly, let us define $K_{n}:X\rightrightarrows X$ by $%
K_{n}(x)=B(K(x);1/(n+n_{0}-1))$ for each $x\in K$ and $n\in \mathbb{N}^{\ast
}$. Since $K$ is almost lower semicontinuous, according to Lemma 2, for each 
$n\in \mathbb{N},$\ there exists a continuous function $f_{n}:X\rightarrow X$%
\ such that $f_{n}(x)\in K_{n}(x)$ for each\textit{\ }$x\in X.$
Brouwer-Schauder fixed point theorem assures that, for each $n\in \mathbb{N},
$ there exists $x_{n}\in X$ such that $x_{n}=f_{n}(x_{n})$ and then, $%
x_{n}\in K_{n}(x_{n}).$

Then, $d(x_{n},K(x_{n}))\rightarrow 0$ when $n\rightarrow \infty $ and since 
$X$ is compact, $\{x_{n}\}$ has a convergent subsequence $\{x_{n_{k}}\}.$
Let $x^{\ast }=\lim_{n_{k}\rightarrow \infty }x_{n_{k}}.$ It follows that $%
d(x^{\ast },K(x_{n_{k}}))\rightarrow 0$ when $n_{k}\rightarrow \infty .$

We claim that $x^{\ast }\in \overline{K}(x^{\ast }).$ Let us assume, to the
contrary, that $x^{\ast }\notin \overline{K}(x^{\ast }).$ $\overline{K}$ has
a closed graph, therefore, its upper and lower sections are closed. We note
that condition i) implies that $\overline{K}$ has nonempty lower sections.
Since $\{x^{\ast }\}\cap \overline{K}^{-1}(x^{\ast })=\emptyset $ and $X$ is
a regular space, there exists $r_{1}>0$ such that $B(x^{\ast };r_{1})\cap 
\overline{K}^{-1}(x^{\ast })=\emptyset $. Consequently, for each $z\in
B(x^{\ast };r_{1}),$ we have that $z\notin \overline{K}^{-1}(x^{\ast }),$
which is equivalent with $x^{\ast }\notin \overline{K}(z)$ or $\{x^{\ast
}\}\cap \overline{K}(z)=\emptyset $. The closedness of each $\overline{K}(z)$
and the regularity of $X$ imply the existence of a real number $r_{2}>0$
such that $B(x^{\ast };r_{2})\cap \overline{K}(z)=\emptyset $ for each $z\in
B(x^{\ast };r_{1}),$ which implies $x_{0}\notin B(\overline{K}(z);r_{2})$
for each $z\in B(x^{\ast };r_{1}).$ Let $r=\min \{r_{1},r_{2}\}.$ Hence, $%
x^{\ast }\notin B(\overline{K}(z);r)$ for each $z\in B(x^{\ast };r),$ and
then, there exists $N^{\ast }\in \mathbb{N}$ such that for each $%
n_{k}>N^{\ast },$ $x^{\ast }\notin B(\overline{K}(x_{n_{k}});r)$ which
contradicts $d(x^{\ast },K(x_{n_{k}}))\rightarrow 0$ as $n\rightarrow \infty 
$. It follows that our assumption is false. Hence, $x^{\ast }\in \overline{K}%
(x^{\ast }).$

Since $T$ has nonempty values, we can pick $y^{\ast }\in T(x^{\ast }).$
Obviously, $G(x^{\ast },y^{\ast })=\emptyset .$ Consequently, the conclusion
holds in Case I.

Case II.

$A=\{(x,y)\in X\times Y:$\textit{\ }$P(x,y)\cap K(x)\neq \emptyset
\}=\{(x,y)\in X\times Y:$\textit{\ }$G(x,y)\neq \emptyset \}$\textit{\ }is
nonempty and open.

Firstly we are proving the convexity of $P(x_{0},y_{0}),$ where $%
(x_{0},y_{0})\in X\times Y$ is arbitrarily fixed$.$ Indeed, let us consider $%
u_{1},u_{2}\in P(x_{0},y_{0})$ and $\lambda \in \lbrack 0,1].$ Since $%
u_{1},u_{2}\in X$ and the set $X$ is convex, the convex combination $%
u=\lambda u_{1}+(1-\lambda )u_{2}$ is an element of $X.$

Further, by using the property of properly $(-C(x_{0}))$ $-$quasi-convexity
of $f(x_{0},y_{0},\cdot )$, we can assume, without loss of generality, that $%
f(x_{0},y_{0},u_{1})\subset f(x_{0},y_{0},u)-C(x_{0}).$

We will prove that $u\in P(x_{0},y_{0}).$ If, to the contrary, $u\notin
P(x_{0},y_{0}),$ then, $f(x_{0},y_{0},u)\nsubseteq -C(x_{0})$ and,
therefore, $f(x_{0},y_{0},u_{1})\subset f(x_{0},y_{0},u)-C(x_{0})\nsubseteq
-C(x_{0})-C(x_{0})\subseteq -C(x_{0}),$ which contradicts $u_{1}\in
P(x_{0},y_{0}).$ Hence, $u\in P(x_{0},y_{0})$ and, consequently, $%
P(x_{0},y_{0})$ is a convex set.

Since $f$ fulfills the property $\mathcal{P}$ with respect to $K$ on $A$,
then, the correspondence $G$ is almost lower semicontinuous on $A$.

Further, let us define the correspondence $H:X\times Y\rightrightarrows X$ by

$H(x,y)=\left\{ 
\begin{array}{c}
G(x,y),\text{ if }(x,y)\in A; \\ 
K(x),\text{ \ \ \ \ \ \ \ otherwise}%
\end{array}%
\right. =\left\{ 
\begin{array}{c}
K(x)\cap P(x,y),\text{ if }(x,y)\in A; \\ 
K(x),\text{ \ \ \ \ \ \ \ \ \ \ \ \ \ \ \ \ \ \ \ otherwise.}%
\end{array}%
\right. $

The correspondence $H$ is nonempty and convex valued.

We notice that for each $(x,y)\in X\times Y,$ $H(x,y)\subseteq K(x)$ and i)
implies that there exists an integer $n_{0}\in \mathbb{N}^{\ast }$ such that 
$B(H(X);\frac{1}{n_{0}})\subset X.$

Now, we are proving that $H$ is almost lower semicontinuous on $X\times Y.$

Let $(x_{0},y_{0})\in X\times Y$ be arbitrary fixed. Let us consider firstly 
$(x_{0},y_{0})\in A.$ Since $H_{\mid A}=G_{\mid A}$, it is nonempty valued.
We conclude that$\ H_{\mid A}$ is almost lower semicontinuous at $%
(x_{0},y_{0})\in A$. $H_{\mid A^{\ast }}=K_{\mid A^{\ast }}$, where $A^{\ast
}=$int$(X\times Y\setminus A),$ and it is nonempty valued. We conclude that$%
\ H_{\mid A^{\ast }}$ is also almost lower semicontinuous on $A^{\ast }$.
Condition iv e) implies that $H$ is almost lower semicontinuous at $%
(x_{0},y_{0})\in $fr$A.$ Hence, $H$ is almost lower semicontinuous on $%
X\times Y.$

For each $n\in \mathbb{N}^{\ast },$ let us define $H_{n}:X\times
Y\rightrightarrows X$ by $H_{n}(x,y)=B(H(x,y);\frac{1}{n+n_{0}-1})$ for each 
$(x,y)\in X\times Y.$

According to Deutsch and Kenderov Lemma, for each $n\in \mathbb{N}^{\ast },$
there exists a continuous function $h_{n}:X\times Y\rightarrow X$ such that $%
h_{n}(x,y)\in H_{n}(x,y)$ for each $(x,y)\in X\times Y.$

For each $n\in \mathbb{N}^{\ast },$ let us define $M_{n}:X\times
Y\rightrightarrows X\times Y$ by $M_{n}(x,y)=(h_{n}(x,y),$cl$T(x))$ for each 
$(x,y)\in X\times Y.$

For each $n\in \mathbb{N}^{\ast },$ $M_{n}$ is closed (since $h_{n}$ is
continuous and cl$T$ is upper semicontinuous with nonempty closed values)
and it has nonempty closed convex values. Then, by the Kakutani-Ky Fan fixed
point Theorem, there exists a fixed point $(x_{n}^{\ast },y_{n}^{\ast })\in
M_{n}(x_{n}^{\ast },y_{n}^{\ast }),$ which implies that $x_{n}^{\ast }\in
H_{n}(x_{n}^{\ast },y_{n}^{\ast })$ and $y_{n}^{\ast }\in $cl$T(x_{n}^{\ast
}).$ We note that $d(x_{n}^{\ast },H(x_{n}^{\ast },y_{n}^{\ast
}))\rightarrow 0$ when $n\rightarrow \infty .$ Since $X\times Y$ is compact, 
$(x_{n}^{\ast },y_{n}^{\ast })_{n\in \mathbb{N}^{\ast }}$ has a convergent
subsequence $(x_{n_{k}}^{\ast },y_{n_{k}}^{\ast })_{n_{k}\in \mathbb{N}%
^{\ast }}.$ Let $(x^{\ast },y^{\ast })=\lim_{n_{k}\rightarrow \infty
}(x_{n_{k}}^{\ast },y_{n_{k}}^{\ast }).$ It follows that $d(x^{\ast
},H(x_{n_{k}}^{\ast },y_{n_{k}}^{\ast }))\rightarrow 0$ when $%
n_{k}\rightarrow \infty .$

The correspondence cl$T$ is closed and therefore, $y^{\ast }\in $cl$%
T(x^{\ast }).$ We claim that $x^{\ast }\in \overline{H}(x^{\ast },y^{\ast }).
$ Let us assume, by the contrary, that $x^{\ast }\notin \overline{H}(x^{\ast
},y^{\ast }),$ or, equivalently, $(x^{\ast },y^{\ast })\notin \overline{H}%
^{-1}(x^{\ast }).$ We note that condition iv f) implies that $\overline{H}$
has nonempty lower sections. In addition, $\overline{H}$ has a closed graph,
therefore, its upper and lower sections are closed. Then, there exists $%
r_{1}>0$ such that $B((x^{\ast },y^{\ast });r_{1})\cap \overline{H}%
^{-1}(x^{\ast })=\emptyset .$ Hence, for each $(x,y)\in B((x^{\ast },y^{\ast
});r_{1}),$ we have that $(x,y)\notin \overline{H}^{-1}(x^{\ast }),$ which
means $x^{\ast }\notin \overline{H}(x,y)$ or $\{x^{\ast }\}\cap \overline{H}%
(x,y)=\emptyset .$ The closedeness of $\overline{H}(x,y)$ implies the
existence of a real number $r_{2}>0$ such that $B(x^{\ast };r_{2})\cap 
\overline{H}(x,y)=\emptyset $ for each $(x,y)\in B((x^{\ast },y^{\ast
});r_{1}),$ which implies $x^{\ast }\notin B(\overline{H}(x,y);r_{2})$ for
each $(x,y)\in B((x^{\ast },y^{\ast });r_{1}).$ Therefore, there exists $%
N^{\ast }\in \mathbb{N}$ such that for each $n_{k}>N^{\ast },$ $x^{\ast
}\notin B(\overline{H}(x_{n_{k}},y_{n_{k}});r_{2}),$ which contradicts $%
d(x^{\ast },H(x_{n}^{\ast },y_{n}^{\ast }))\rightarrow 0$ when $%
n_{k}\rightarrow \infty .$ We conclude that our assumption is false. Hence, $%
x^{\ast }\in \overline{H}(x^{\ast },y^{\ast }).$

The next step of the proof is to show that $(x^{\ast },y^{\ast })\notin A.$
We firstly note that, according to $iv$ $a)$, $x\notin \overline{P}(x,y)$
for each $(x,y)\in X\times Y.$ Suppose, to the contrary, that $(x^{\ast
},y^{\ast })\in A.$ In this case, $x^{\ast }\in \overline{H}(x^{\ast
},y^{\ast })=\overline{G}(x^{\ast },y^{\ast })\subseteq \overline{P}(x^{\ast
},y^{\ast }),$ which contradicts the assertion above. Thus, $(x^{\ast
},y^{\ast })\notin A.$ Therefore, $x^{\ast }\in \overline{K}(x^{\ast }),$ $%
y^{\ast }\in $cl$T(x^{\ast })$ and $G(x^{\ast },y^{\ast })=K(x^{\ast })\cap
P(x^{\ast },y^{\ast })=\emptyset .$

Consequently, there exist $x^{\ast },y^{\ast }\in X$ such that $x^{\ast }\in 
\overline{K}(x^{\ast }),$ $y^{\ast }\in $cl$T(x^{\ast })$ and $f(x^{\ast
},y^{\ast },u)\nsubseteq -C(x^{\ast })$ for each $u\in K(x^{\ast }).$

We note that $f(x^{\ast },y^{\ast },u)\nsubseteq -C(x^{\ast })$ implies $%
f(x^{\ast },y^{\ast },u)\nsubseteq -$int$C(x^{\ast })$ and then, the pair $%
(x^{\ast },y^{\ast })\in X\times Y$\textit{\ }is a solution of
GVQEP.\smallskip

Now, we recall the following important results.

\begin{definition}
Let $K$ be a subset of a topological vector space $E.$ According to Had\v{z}i%
\'{c} (\cite{ha}), $K$ is said to be of the Zima type if
\end{definition}

(Z) \textit{for each neighborhood }$U$ of $0$ in $E,$ \textit{there exists a
neighborhood }$V$\textit{\ of the origin in }$E$\textit{\ such that} $%
co(V\cap (E-E))\subset U.$\smallskip

In \cite{park2004}, Park established the following result.

\begin{theorem}
(Theorem 3.1 in \cite{park2004}) Let $X$ be a convex subset of a topological
vector space $E.$ Let $T:X\rightrightarrows X$ be an upper semicontinuous
(respectively a lower semicontinuous) correspondence with nonempty convex
values such that $T(X)$ is of the Zima type.
\end{theorem}

\textit{If }$T(X)$\textit{\ is totally bounded, then, for each neighborhood }%
$U$\textit{\ of }$0$\textit{\ in }$E$\textit{, there exists a }$U-$\textit{%
almost fixed point of }$T,$ \textit{that is, a point }$x_{U}\in X$\textit{\
such that }$T(x_{U})\cap (x_{U}+U)\neq \emptyset .\medskip $

Theorem 3 is stated in terms of upper semicontinuity for $f(\cdot ,\cdot
,u):X\times Y\rightrightarrows Z$ $(u\in X)$ and lower semicontinuity for $K$
and $T.$ By using the Park Theorem, we firstly prove the existence of the $%
U- $almost fixed points for a correspondence we construct. An approximation
technique is used in order to prove our statement. We work on regular vector
topological spaces (we drop the locally convexity). We also drop the
property of properly $C-$ quasi convexity of $f(x,y,\cdot ).$ All of these
make the differences between Theorem 3 and Theorem 3.2.4 in \cite{lin3}. In
fact, our result is an improvement of the mentioned theorem. We underline
that we use a different argument for the proof.

\begin{theorem}
Let $Z$ be a Hausdorff topological vector space and l\textit{et }$X,$ $Y$%
\textit{\ be nonempty compact convex subsets of regular topological vector
spaces }$E$\textit{, respectively }$F$\textit{. } \textit{Assume that:}
\end{theorem}

\textit{i) }$K:X\rightrightarrows X$\textit{\ is a lower semicontinuous
correspondence with nonempty convex open values such that }$K(X)$ \textit{is
of the Zima type and }$\overline{K}$ \textit{has nonempty lower sections; }

\textit{ii) }$T:X\rightrightarrows Y$\textit{\ is a lower semicontinuous
correspondence with nonempty convex values such that }$T(X)$ \textit{is of
the Zima type and }$\overline{T}$\textit{\ has nonempty lower sections; }

\textit{iii) }$C:X\rightrightarrows Z$\textit{\ is a correspondence such
that, for all }$x\in X,$\textit{\ }$C(x)$\textit{\ is a pointed closed
convex cone with int}$C(x)\neq \emptyset ;$

\textit{iv) }$f:X\times Y\times X\rightrightarrows Z$\textit{\ is a
correspondence with nonempty values which satisfies the following conditions:%
}

\textit{\ \ \ \ \ \ \ \ \ \ a) for all }$(x,y)\in X\times Y,$\textit{\ }$%
(x,y,x)\notin $cl$\{(x,y,z)\in X\times Y\times X:z\in $co$\{K(x)\cap \{u\in
X:f(x,y,u)\subseteq -$int$C(x)\}\}\};$

\textit{\ \ \ \ \ \ \ \ b) }$A=\{(x,y)\in X\times Y:$\textit{\ there exists }%
$u\in K(x)$ \textit{such that} \textit{\ }$f(x,y,u)\subseteq -$\textit{int}$%
C(x)\}$\textit{\ is closed or empty;}

\ \ \ \ \ \ \ \ \ \textit{c)} \textit{for each }$u\in X,$ $f(\cdot ,\cdot
,u):X\times Y\rightrightarrows Z$\textit{\ is upper semicontinuous on }$A$%
\textit{;}

\ \ \ \ \ \ \ \ \textit{d) the correspondence }$Q:X\rightrightarrows Z$%
\textit{\ defined by }$Q(x)=Z\backslash (-C(x))$\textit{\ for each }$x\in X$ 
\textit{is upper semicontinuous on pr}$_{X}A$ \textit{(if} $A\neq \emptyset 
\mathit{)}$\textit{;}

\ \ \ \ \ \ \ \ \textit{e)} \textit{for each }$z\in X,$ $\{(x,y)\in X\times
Y:(x,y,z)\in $cl$_{X\times Y\times X}\{(x,y,u)\in X\times Y\times
X:f(x,y,u)\subseteq -$int$C(x)$ \textit{and} $u\in K(x)\}\}$ \textit{is
nonempty.}

\textit{Then, there exists a solution }$(x^{\ast },y^{\ast })\in X\times Y$%
\textit{\ of GVQEP.\medskip }

\textit{Proof.} The key idea of the proof is to use the Park Theorem.
Towards this end, we need to construct the following correspondences.

Let $P,G:X\times Y\rightrightarrows X$ be defined by

$P(x,y)=\{u\in X:f(x,y,u)\subseteq -$int$C(x)\},$ for each $(x,y)\in X\times
Y$ and

$G(x,y)=K(x)\cap P(x,y),$ for each $(x,y)\in X\times Y.$

We are going to show that there exists a pair\textit{\ }$(x^{\ast },y^{\ast
})\in X\times Y$ such that $x^{\ast }\in \overline{K}(x^{\ast }),$ $y^{\ast
}\in \overline{T}(x^{\ast })$ and $K(x^{\ast })\cap P(x^{\ast },y^{\ast
})=\emptyset .$

We consider two cases.

Case I.

$A=\{(x,y)\in X\times Y:$\textit{\ }$P(x,y)\cap K(x)\neq \emptyset
\}=\{(x,y)\in X\times Y:$\textit{\ }$G(x,y)\neq \emptyset \}\mathit{=}%
\emptyset $.

In this case, $G(x,y)=\emptyset $ for each $(x,y)\in X\times Y.$ We can
apply the Park fixed point Theorem to $K.$

For each symmetric open neighborhood $W$\ of $0$\ in $E$, there exists a
symmetric open neighborhood $U$\ of $0$\ in $E$ such that $U+U\subset W.$
According to the Park Theorem, for each such a neighborhood $U,$ there
exists points $x_{U},$ $y_{U}\in X$\ such that\textit{\ }$x_{U}\in
T(x_{U})+U $ and $x_{U}\in y_{U}+U.$

Since $X$ is compact, $\{x_{U}\}$ has a convergent subsequence $%
\{x_{U^{\prime }}\}.$ Let $x^{\ast }$ be the limit of $\{x_{U^{\prime }}\}.$
It follows that $x^{\ast }\in T(x_{U^{\prime }})+W^{\prime }$ for each
symmetric open neighborhood $U^{\prime }$ of $0$ with the property that $%
U^{\prime }+U^{\prime }\subset W^{\prime }.$

Let us assume that $x^{\ast }\notin \overline{K}(x^{\ast }).$We note that
condition i) implies that $\overline{K}$ has nonempty lower sections. In
addition, $\overline{K}$ has a closed graph, therefore, its upper and lower
sections are closed. Since $\{x^{\ast }\}\cap \overline{K}^{-1}(x^{\ast
})=\emptyset $ and $X$ is a regular space, there exists $V_{1}$ an open
neighborhood of $0$ such that $(x^{\ast }+V_{1})\cap \overline{K}%
^{-1}(x^{\ast })=\emptyset $. Consequently, for each $z\in (x^{\ast
}+V_{1}), $ we have that $z\notin \overline{K}^{-1}(x^{\ast }),$ which is
equivalent with $x^{\ast }\notin \overline{K}(z)$ or $\{x^{\ast }\}\cap 
\overline{K}(z)=\emptyset $. The closedness of each $\overline{K}(z)$ and
the regularity of $X$ imply the existence of $V_{2},$ an open neighborhood
of $0,$ such that $(x^{\ast }+V_{2})\cap \overline{K}(z)=\emptyset $ for
each $z\in x^{\ast }+V_{1},$ which implies $x^{\ast }\notin \overline{K}%
(z)+V_{2}$ for each $z\in x^{\ast }+V_{1}.$ Let $V=V_{1}\cap V_{2}.$ Hence, $%
x^{\ast }\notin \overline{K}(z)+V$ for each $z\in x^{\ast }+V,$ and then,
there exists $U^{\ast },$ an open neighborhood of 0 such that for each
symmetric open neighborhood of $0,$ $U^{\prime },$ with the property that $%
U^{\prime }\subset U^{\ast },$ it is true that $x^{\ast }\notin \overline{K}%
(x_{U^{\prime }})+V$ and therefore, $x^{\ast }\notin \overline{K}%
(x_{U^{\prime }\cap V})+W^{\prime }\cap V.$ The last assertion contradicts $%
x^{\ast }\in K(x_{U^{\prime }})+W^{\prime }$ for each symmetric open
neighborhood $U^{\prime }$ of $0$ with the property that $U^{\prime
}+U^{\prime }\subset W^{\prime }.$ It follows that our assumption is false.
Hence, $x^{\ast }\in \overline{K}(x^{\ast }).$

Since $T$ has nonempty values, we can pick $y^{\ast }\in T(x^{\ast }).$
Obviously, $G(x^{\ast },y^{\ast })=\emptyset .$ Consequently, the conclusion
holds in Case I.

Case II.

$A=\{(x,y)\in X\times Y:$\textit{\ }$P(x,y)\cap K(x)\neq \emptyset
\}=\{(x,y)\in X\times Y:$\textit{\ }$G(x,y)\neq \emptyset \}$\textit{\ }is
nonempty and closed.

In this case, we firstly show that $P$ is lower semicontinuous on $A$ with
nonempty values$.$ In order to do this, we are going to show that for each
closed set $V$ in $X,$ $W_{0}=\{(x,y)\in A:P(x,y)\subset V\}$ is closed.

We notice that $W_{0}=\{(x,y)\in A:P(x,y)\subset V\}=\{(x,y)\in A:\{u\in
X:f(x,y,u)\subseteq -$int$C(x)\}\subset V\}.$

Let $(x_{n},y_{n})_{n\in \mathbb{N}}\subset W_{0},$ such that $%
(x_{n},y_{n})\rightarrow (x_{0},y_{0}).$ It implies that for each $n\in 
\mathbb{N},$ $\{u\in X:f(x_{n},y_{n},u)\subseteq -$int$C(x_{n})\}\subset V.$
We want to prove that $(x_{0},y_{0})\in W_{0},$ which is equivalent with $%
\{u\in X:f(x_{0},y_{0},u)\subseteq -$int$C(x_{0})\}\subset V.$ If there
exists $n_{0}\in \mathbb{N}$ such that $\{u\in X:f(x_{0},y_{0},u)\subseteq -$%
int$C(x_{0})\}\subset \{u\in X:f(x_{n_{0}},y_{n_{0}},u)\subseteq -$int$%
C(x_{n_{0}})\},$ the last assertion is true.

Let $u_{0}\in X,$ with the property that $f(x_{0},y_{0},u_{0})\subseteq -$int%
$C(x_{0}).$ Since $f(\cdot ,\cdot ,u_{0})$ is upper semicontinuous on $A$,
there exists an open neighborhood $W=U\times V$ of $(x_{0},y_{0})$ in $A$
such that $f(x,y,u_{0})\subseteq -$int$C(x_{0})$ for each $(x,y)\in U\times
V.$ Then, there exists $N_{1}\in \mathbb{N}$ such that $f(x_{n},y_{n},u_{0})%
\subseteq -$int$C(x_{0})$ for each $n>N_{1}.$ The correspondence $%
Q:X\rightrightarrows Z$ is upper semicontinuous on pr$_{X}A$ and then, there
exists an open neighborhood $U_{0}$ of $x_{0}$ in pr$_{X}A$ such that for
each $x\in U_{0},$ $Q(x)\subseteq Q(x_{0})$ or, $-C(x_{0})\subseteq -C(x)$.
Then, there exists $N_{2}\in \mathbb{N}$ such that $-C(x_{0})\subseteq
-C(x_{n})$ for each $n>N_{2}.$ Let $N=\max (N_{1}.N_{2}).$ Therefore, $%
f(x_{n},y_{n},u_{0})\subseteq -$int$C(x_{n})$ for each $n>N.$ Consequently, $%
\{u\in X:f(x_{0},y_{0},u)\subseteq -$int$C(x_{0})\}\subset \{u\in
X:f(x_{n},y_{n},u)\subseteq -$int$C(x_{n})\}$ for each $n>N$ and then, $%
\{u\in X:f(x_{0},y_{0},u)\subseteq -$int$C(x_{0})\}\subset V.$ We conclude
that $(x_{0},y_{0})\in W_{0},$ $W_{0}$ is closed and $P$ is lower
semicontinuous on $A.$

The correspondence $G:A\rightrightarrows X,$ defined by $G(x,y)=K(x)\cap
P(x,y)$ for each $(x,y)\in A,$ is lower semicontinuous on $A,$ since $K$ and 
$P$ are lower semicontinuous with nonempty values, $K$ has open values and $%
G $ is nonempty valued.

Further, let us define the correspondence $H:X\times Y\rightrightarrows X$ by

$H(x,y)=\left\{ 
\begin{array}{c}
\text{co}G(x,y),\text{ if }(x,y)\in A; \\ 
K(x),\text{ \ \ \ \ \ \ \ otherwise}%
\end{array}%
\right. =\left\{ 
\begin{array}{c}
\text{co}(K(x)\cap P(x,y)),\text{ if }(x,y)\in A; \\ 
K(x),\text{ \ \ \ \ \ \ \ \ \ \ \ \ \ \ \ \ \ \ \ otherwise.}%
\end{array}%
\right. $

According to Lemma 1, the correspondence $H$ is lower semicontinuous on $%
X\times Y$. In addition, it has nonempty and convex values.

If $K(X)$ is of the Zima type and $H(x,y)\subseteq K(x)$ for each $(x,y)\in
X\times Y,$ it can be shown easily that $H(X\times Y)$ is also of the Zima
type.

Let us define $M:X\times Y\rightrightarrows X\times Y$ by $%
M(x,y)=(H(x,y),T(x))$ for each $(x,y)\in X\times Y.$

The correspondence $M$ is lower semicontinuous with nonempty convex values.
We note that $M(X,Y)$ is totally bounded.

We will prove that $M$ is of the Zima type$.$ We notice that $M(X\times
Y)=H(X\times Y)\times T(X).$ We know that the sets $H(X\times Y)$ and $T(X)$
are of the Zima type and therefore, for each neighborhoods $U_{1}$ and $%
U_{2} $ of $0$ in $E$, respectively $F$, there exist the neighborhoods $%
V_{1} $ and $V_{2\text{ }}$ of $0$ in $E$, respectively $F$, such that

co$(V_{1}\cap (H(X\times Y)-H(X\times Y))\subset U_{1}$ and co$(V_{2}\cap
(T(X)-T(X))\subset U_{2}.$

Therefore, co$((V_{1}\times V_{2})\cap (H(X\times Y)\times T(X)-H(X\times
Y)\times T(X))\subseteq $

$\subseteq $co$((V_{1}\times V_{2})\cap (H(X\times Y)-H(X\times Y)\times
(T(X)-T(X)))\subseteq $

$\subseteq $co$(V_{1}\cap (H(X\times Y)-H(X\times Y))\times (V_{2}\cap
(T(X)-T(X))\subseteq $

$\subseteq $co$(V_{1}\cap (H(X\times Y)-H(X\times Y))\times $co$(V_{2}\cap
(T(X)-T(X))\subseteq U_{1}\times U_{2}.$

It follows that for each for each neighborhood $U_{1}\times U_{2}$ of $(0,0)$
in $E\times F$, there exists $a$ neighborhood $V_{1}\times V_{2\text{ }}$ of 
$(0,0)$ in $E\times F$, such that co$((V_{1}\times V_{2})\cap (H(X\times
Y)\times T(X)-H(X\times Y)\times T(X))\subseteq U_{1}\times U_{2},$ that is, 
$M$ is of the Zima type$.$

Then, according to the Park Theorem, there exists a $U$-almost fixed point \
of $M,$ that is, for each neighborhood $W=(U,V)$ of $(0,0)$ in $X\times Y,$
there exists a point $(x_{U},y_{V})$ such that $M(x_{U},y_{V})\cap
((x_{U},y_{V})+(U,V))\neq \emptyset .$ It follows that $H(x_{U},y_{V})\cap
(x_{U}+U)\neq \emptyset $ and $T(x_{U})\cap (y_{V}+V)\neq \emptyset .$

For each symmetric open neighborhood $W=(W_{1},W_{2})$ of $(0,0)$\ in $%
E\times F$, there exist symmetric open neighborhoods $U$ and $V$\ of \ the
origin\ in $E,$ respectively $F,$ such that $(U,V)+(U,V)\subset W.$ Then$,$
there exist points $x_{U},$ $y_{V},z_{V}\in X$\ such that\textit{\ }$%
z_{V}\in T(x_{U})+V$ and $z_{V}\in y_{V}+V.$

$X$ and $Y$ are compact sets. Therefore, $\{x_{U}\}_{U}$ and $\{y_{V}\}_{V}$
have convergent subsequences $\{x_{U^{\prime }}\},$ respectively, $%
\{y_{V^{\prime }}\}.$ Let $x^{\ast }$ be the limit of $\{x_{U^{\prime }}\}$
and $y^{\ast }$ be the limit of $\{y_{V^{\prime }}\}.$ It follows that $%
y^{\ast }\in T(x_{U^{\prime }})+W_{2}^{\prime }$ for each symmetric open
neighborhood $U^{\prime }$ of $0$ in $E$ with the property that $(U^{\prime
},V^{\prime })+(U^{\prime },V^{\prime })\subset W^{\prime },$ where $%
V^{\prime }$ is a symmetric open neighborhood of $0$ in $F$ and $W^{\prime
}=(W_{1}^{\prime },W_{2}^{\prime })$ is a symmetric open neighborhood of $%
(0,0)$\ in $E\times F.$ We claim that $y^{\ast }\in \overline{T}(x^{\ast }).$

\QTP{Body Math}
Let us assume, to the contrary, that $y^{\ast }\notin \overline{T}(x^{\ast
}).$ We note that condition ii) implies that $\overline{T}$ has nonempty
lower sections. In addition, $\overline{T}$ has a closed graph, therefore,
its upper and lower sections are closed. Since $\{x^{\ast }\}\cap \overline{T%
}^{-1}(y^{\ast })=\emptyset $ and $E$ is a regular space, there exists $%
U_{0} $ an open neighborhood of $0$ in $E$ such that $(x^{\ast }+U_{0})\cap 
\overline{T}^{-1}(y^{\ast })=\emptyset $. In this case, we remark that, for
each $z\in (x^{\ast }+U_{0}),$ we have $z\notin \overline{T}^{-1}(y^{\ast
}). $ As result, we notice that $y^{\ast }\notin \overline{T}(z)$ or $%
\{y^{\ast }\}\cap \overline{T}(z)=\emptyset $. The closedness of each $%
\overline{T}(z)$ and the regularity of $F$ imply the existence of $V_{0},$
an open neighborhood of $0$ in $F,$ such that $(y^{\ast }+V_{0})\cap 
\overline{T}(z)=\emptyset $ for each $z\in x^{\ast }+U_{0},$ which further
implies $y^{\ast }\notin \overline{T}(z)+V_{0}$ for each $z\in x^{\ast
}+U_{0}.$ Then, there exists $U^{\ast },$ an open neighborhoods of the
origin in $E,$ such that for each symmetric open neighborhood $U^{\prime }$
of \ the origin in $E,$ with the property that $U^{\prime }\subset U^{\ast
}, $ it is true that $y^{\ast }\notin \overline{T}(x_{U^{\prime }})+V_{0}$
and therefore, $y^{\ast }\notin \overline{T}(x_{U^{\prime }\cap
U_{0}})+W_{2}^{\prime }\cap V_{0}.$ The last assertion contradicts $y^{\ast
}\in \overline{T}(x_{U^{\prime }})+W_{2}^{\prime }$ for each symmetric open
neighborhood $U^{\prime }$ of $0$ in $E$ with the property that $(U^{\prime
},V^{\prime })+(U^{\prime },V^{\prime })\subset W^{\prime },$ where $%
V^{\prime }$ is a symmetric open neighborhood of $0$ in $F$ and $W^{\prime
}=(W_{1}^{\prime },W_{2}^{\prime })$ is a symmetric open neighborhood of $%
(0,0)$\ in $E\times F.$ It follows that our assumption is false. Hence, $%
y^{\ast }\in \overline{T}(x^{\ast }).$

Similarly, we can show that $x^{\ast }\in \overline{H}(x^{\ast },y^{\ast }).$
For this, we take into account that $iv$ $f)$ implies that $\overline{H}$
has nonempty lower sections.

The next step of the proof is to show that $(x^{\ast },y^{\ast })\notin A.$
Suppose, to the contrary, that $(x^{\ast },y^{\ast })\in A.$ Then, $x^{\ast
}\in \overline{H}(x^{\ast },y^{\ast })=\overline{G}(x^{\ast },y^{\ast }),$
which contradicts $iv$ $a).$ Thus, $(x^{\ast },y^{\ast })\notin A.$
Therefore, $x^{\ast }\in \overline{K}(x^{\ast }),$ $y^{\ast }\in \overline{T}%
(x^{\ast })$ and $G(x^{\ast },y^{\ast })=K(x^{\ast })\cap P(x^{\ast
},y^{\ast })=\emptyset .$

Consequently, there exist $x^{\ast },y^{\ast }\in X$ such that $x^{\ast }\in 
\overline{K}(x^{\ast }),$ $y^{\ast }\in \overline{T}(x^{\ast })$ and $%
f(x^{\ast },y^{\ast },u)\nsubseteq -$int$C(x^{\ast })$ for each $u\in
K(x^{\ast }).$

\begin{remark}
Instead of assuming that $K$ has open values, we can assume that for each $%
(x,y)\in A,$ $\{u\in X:f(x,y,u)\subseteq -$int$C(x)\}$ is open, that is, $%
P_{\mid A}$ has open values. One of these two conditions is needed in order
to obtain the lower semicontinuity of $P_{\mid A}.$
\end{remark}

\section{Applications of the KKM principle to the equilibrium problems}

Many essential results of the equilibrium theory can be derived from the KKM
principle. We recall the KKM principle here. We note that its open version
is due to Kim \cite{kim1987} and Shih and Tan \cite{shih}.\medskip

Let $X$ be a subset of a topological vector space and $D$ a nonempty subset
of $X$ such that co$D\subset X.$

$T:D\rightrightarrows X$ is called a \textit{KKM correspondence} if co$%
N\subset T(N)$ for each $N\in \langle D\rangle ,$ where $\langle D\rangle $
denotes the class of all nonempty finite subsets of $D$.\medskip

\textbf{KKM principle }Let $D$ be a set of vertices of a simplex $S$ and $%
T:D\rightrightarrows 2^{S}$ a correspondence with closed (respectively open)
values such that

co$N\subset T(N)$ for each $N\subset D.$

Then, $\tbigcap\nolimits_{z\in D}T(z)\neq \emptyset .\medskip $

The following lemma is a consequence of the KKM principle. It will be used
to obtain new existence results for generalized vector equilibrium problems
in this section.

\begin{lemma}
Let $X$ be a subset of a topological vector space, $D$ a nonempty subset of $%
X$ such that co$D\subset X$ and $T:D\rightrightarrows X$ a \textit{KKM
correspondence with closed (respectively open) values. Then }$\{T(z)\}_{z\in
D}$ has the finite intersection property.\smallskip
\end{lemma}

In this section, we approach a particular case of the problem considered in
Section 3.\medskip

Let $X$ be a nonempty subset of a topological vector space $E$ and let $Z$
be a topological vector space. Let $C\subset Z$\textit{\ }be a pointed
closed convex cone with nonempty interior$.$ Let $K:X\rightrightarrows X$
and $f:X\times X\rightrightarrows Z$ be correspondences with nonempty
values. We consider the following generalized vector quasi-equilibrium
problem:

(GVQEP 2): Find $x^{\ast }\in X$ such that $x^{\ast }\in K(x^{\ast })$ and $%
f(x^{\ast },y)\nsubseteq -$int$C$ for each $y\in K(x^{\ast })$.$\medskip $

Now, we are establishing an existence theorem for a generalized vector
equilibrium problem by using Lemma 3.

\begin{theorem}
\textit{Let }$Z$\textit{\ be a Hausdorff topological vector space and let }$%
X $\textit{\ be a nonempty compact convex subset of a\ topological vector
space }$E$\textit{. Let }$C\subset Z$\textit{\ be a pointed closed convex
cone with nonempty interior and let }$K:X\rightrightarrows X$\textit{\ and }$%
f:X\times X\rightrightarrows Z$\textit{\ be correspondences with nonempty
values. Assume that:}
\end{theorem}

\textit{i) }$K$\textit{\ is open valued and for each }$x\in X,$\textit{\ the
set }$\{u\in X:f(x,u)\subseteq -$\textit{int}$C\}$\textit{\ is open; }

\textit{ii) \ }$f(x,x)\nsubseteq -$\textit{int}$C$\textit{\ for each }$x\in
X $\textit{; }

\textit{iii) there exists }$M\in \langle X\rangle $\textit{\ such that }$%
\bigcup\nolimits_{x\in M\cap A}[\{u\in X:f(x,u)\subseteq -$\textit{int}$%
C\}\cap K(x)]\bigcup \bigcup\nolimits_{x\in M\backslash A}K(x)=X,$\textit{\
where }$A=\{x\in X:$\textit{\ there exists }$u\in K(x)$\textit{\ such that }$%
f(x,u)\subseteq -$\textit{int}$C\}$\textit{;}

\textit{iv) }$K^{-1}:X\rightrightarrows X$\textit{\ is convex valued and for
each }$u\in X,$\textit{\ the set }$\{x\in X:f(x,u)\subseteq -$\textit{int}$%
C\}\cup (X\backslash A)$\textit{\ is convex}$;$

\textit{Then, there exists }$x^{\ast }\in X$\textit{\ such that }$x^{\ast
}\in K(x^{\ast })$\textit{\ and }$f(x^{\ast },y)\nsubseteq -$\textit{int}$C$%
\textit{\ for each }$y\in K(x^{\ast })$\textit{.\medskip }

\textit{Proof.} Let $P,G:X\rightrightarrows X$ be defined by

$P(x)=\{u\in X:f(x,u)\subseteq -$int$C\},$ for each $x\in X$ and

$G(x)=K(x)\cap P(x),$ for each $x\in X.$

We are going to show that there exists $x^{\ast }\in X$ such that $x^{\ast
}\in K(x^{\ast })$ and $K(x^{\ast })\cap P(x^{\ast })=\emptyset .$

We consider two cases.

Case I.

$A=\{x\in X:$\textit{\ }$P(x)\cap K(x)\neq \emptyset \}=\{x\in X:$\textit{\ }%
$G(x)\neq \emptyset \}$\textit{\ }is nonempty.

The correspondence $G:A\rightrightarrows X,$ defined by $G(x)=K(x)\cap P(x)$
for each $x\in A,$ is nonempty valued on $A.$ We note that for each $u\in X,$
$G^{-1}(u)=P^{-1}(u)\cap K^{-1}(u)$ is a convex set as intersection of
convex sets.

Further, let us define the correspondence $H:X\rightrightarrows X$ by

$H(x)=\left\{ 
\begin{array}{c}
G(x),\text{ if }x\in A; \\ 
K(x),\text{\ otherwise}%
\end{array}%
\right. =\left\{ 
\begin{array}{c}
K(x)\cap P(x),\text{ if }x\in A; \\ 
K(x),\text{ \ \ \ \ \ \ \ \ \ otherwise.}%
\end{array}%
\right. $

According to i), $H$ is open valued and according to iv), for each $u\in X,$ 
$P^{-1}(u)\cup (X\backslash A)$ is convex.

For each $u\in X,$

$H^{-1}(u)=\{x\in X:$\textit{\ }$\mathit{u}\in H(x)\}=$

\ \ \ \ \ \ \ \ \ \ \ $=\{x\in A:$\textit{\ }$\mathit{u}\in G(x)\}\cup
\{x\in X\setminus A:$\textit{\ }$\mathit{u}\in K(x)\}=$

\ \ \ \ \ \ \ \ \ \ \ $=\{x\in A:$\textit{\ }$\mathit{u}\in K(x)\cap
P(x)\}\cup \{x\in X\setminus A:\mathit{\ u}\in K(x)\}$

$\ \ \ \ \ \ \ \ =[A\cap P^{-1}(u)\cap K^{-1}(u)]\cup \lbrack (X\backslash
A)\cap $\textit{\ }$K^{-1}(u)]=$

\ \ \ \ \ \ \ \ \ $=[P^{-1}(u)\cup (X\backslash A)]\cap K^{-1}(u).$

$H^{-1}(u)$ is convex, since $P^{-1}(u)\cup (X\backslash A)$ and $K^{-1}(u)$
are convex sets.

Assumption iii) implies that there exists $M\in \langle X\rangle $ such that 
$\tbigcup\nolimits_{x\in M}H(x)=X.$

Let us define $F:X\rightrightarrows X$ by $F(x):=X\backslash H(x)$ for each%
\textit{\ }$x\in X.$

Then, $F$ is closed valued and $\tbigcap\nolimits_{x\in M}F(x)=X\backslash
\tbigcup\nolimits_{x\in M}(H(x))=\emptyset .$

According to Lemma 3, we conclude that $F$ is not a KKM correspondence.
Thus, there exists $N\in \langle X\rangle $ such that\textit{\ }co$%
N\varsubsetneq F(N)=\tbigcup\nolimits_{x\in N}(X\backslash H(x)).$

Hence, there exists $x^{\ast }\in $co$N$ with the property that $x^{\ast
}\in H(x)$ for each $x\in N.$ Therefore, there exists $x^{\ast }\in $co$N$
such that $x^{\ast }\in H(x)$ for each $x\in N,$ which implies $N\subset
H^{-1}(x^{\ast })$. Further, it is true that co$N\subset $co$H^{-1}(x^{\ast
})=H^{-1}(x^{\ast })$. Consequently, $x^{\ast }\in $co$N\subset $co$%
H^{-1}(x^{\ast })=H^{-1}(x^{\ast }),$ which means that $x^{\ast }\in
H(x^{\ast }),$ that is, $x^{\ast }$ is a fixed point for $H.$

We notice that, according to $ii)$, $x\notin P(x)$ for each $x\in X,$ and
then, $x^{\ast }\notin A.$ Therefore, $x^{\ast }\in K(x^{\ast })$ and $%
G(x^{\ast })=K(x^{\ast })\cap P(x^{\ast })=\emptyset .$

Consequently, there exists $x^{\ast }\in X$ such that $x^{\ast }\in
K(x^{\ast })$ and $f(x^{\ast },u)\nsubseteq -$int$C$ for each $u\in
K(x^{\ast }).$

Case II.

$A=\{x\in X:$\textit{\ }$P(x)\cap K(x)\neq \emptyset \}=\{x\in X:$\textit{\ }%
$G(x)\neq \emptyset \}\mathit{=}\emptyset $.

In this case, $G(x)=\emptyset $ for each $x\in X.$

Let us define $F:X\rightrightarrows X$ by $F(x):=X\backslash K(x)$ for each%
\textit{\ }$x\in X.$ The proof follows the same line as above and we obtain
that there exists $x^{\ast }\in X$ such that $x^{\ast }\in K(x^{\ast }).$

Obviously, $G(x^{\ast })=\emptyset .$ Consequently, the conclusion holds in
Case I.

\begin{remark}
Assumption i) can be replaced with
\end{remark}

i') $K$\ is closed valued and for each $x\in X,$ the set $\{u\in
X:f(x,u)\subseteq -$int$C\}$ is closed.\textit{\ }

In this case, $F$ is open valued.\smallskip

The next theorem can be obtained easily by following a similar reasoning as
in the above result.

\begin{theorem}
\textit{Let }$Z$\textit{\ be a Hausdorff topological vector space and let }$%
X $\textit{\ be a nonempty compact convex subset of a\ topological vector
space }$E$\textit{. Let }$C\subset Z$\textit{\ be a pointed closed convex
cone with nonempty interior and let }$K:X\rightrightarrows X$\textit{\ and }$%
f:X\times X\rightrightarrows Z$\textit{\ be correspondences with nonempty
values. Assume that:}
\end{theorem}

\textit{i) }$K$\textit{\ is open (resp. closed) valued and for each }$x\in
X, $\textit{\ the set }$\{u\in X:f(x,u)\cap $\textit{int}$C=\emptyset \}$%
\textit{\ is open (resp. closed); }

\textit{ii) \ }$f(x,x)\cap $\textit{int}$C=\emptyset $\textit{\ for each }$%
x\in X$\textit{; }

\textit{iii) there exists }$M\in \langle X\rangle $\textit{\ such that }$%
\bigcup\nolimits_{x\in M\cap A}[\{u\in X:f(x,u)\cap $\textit{int}$%
C=\emptyset \}\cap K(x)]\bigcup \bigcup\nolimits_{x\in M\backslash A}K(x)=X,$%
\textit{\ where }$A=\{x\in X:$\textit{\ there exists }$u\in K(x)$\textit{\
such that }$f(x,u)\cap $\textit{int}$C=\emptyset \}$\textit{;}

\textit{iv) }$K^{-1}:X\rightrightarrows X$\textit{\ is convex valued and for
each }$u\in X,$\textit{\ the set }$\{x\in X:f(x,u)\cap $\textit{int}$%
C=\emptyset \}\cup (X\backslash A)$\textit{\ is convex}$;$

\textit{Then, there exists }$x^{\ast }\in X$\textit{\ such that }$x^{\ast
}\in K(x^{\ast })$\textit{\ and }$f(x^{\ast },y)\cap $\textit{int}$C\neq
\emptyset $\textit{\ for each }$y\in K(x^{\ast })$\textit{.\medskip }

\textit{Proof.} Let $P,G:X\rightrightarrows X$ be defined by

$P(x)=\{u\in X:f(x,u)\cap $int$C=\emptyset \},$ for each $x\in X$ and

$G(x)=K(x)\cap P(x),$ for each $x\in X.$

We are going to show that there exists $x^{\ast }\in X$ such that $x^{\ast
}\in K(x^{\ast })$ and $K(x^{\ast })\cap P(x^{\ast })=\emptyset .$

The rest of the proof follows the same line as the proof of Theorem
4.\smallskip

We obtain a new theorem concerning the existence of solutions for a
generalized vector equilibrium problem.

\begin{theorem}
\textit{Let }$Z$\textit{\ be a Hausdorff topological vector space and let }$%
X $\textit{\ be a nonempty compact convex subset of a\ topological vector
space }$E$\textit{. Let }$C\subset Z$\textit{\ be a pointed closed convex
cone with nonempty interior and let }$K:X\rightrightarrows X$\textit{\ and }$%
f:X\times X\rightrightarrows Z$\textit{\ be correspondences with nonempty
values. Assume that:}
\end{theorem}

\textit{i) }$K$\textit{\ is open valued and for each }$x\in X,$\textit{\ the
set }$\{u\in X:f(x,u)\subseteq -$\textit{int}$C\}$\textit{\ is open; }

\textit{ii) \ }$f(x,x)\nsubseteq -$\textit{int}$C$\textit{\ for each }$x\in
X $\textit{; }

\textit{iii) if }$A=\{x\in X:$\textit{\ there exists }$u\in K(x)$\textit{\
such that }$f(x,u)\subseteq -$\textit{int}$C\},$\textit{\ then for each }$%
N\in \langle A\rangle ,$\textit{\ }$($\textit{co}$N\smallsetminus N)\cap
A=\emptyset ;$

\textit{iv) there exists }$M\in \langle A\rangle $\textit{\ such that }$%
\bigcup\nolimits_{x\in M}\{u\in X:f(x,u)\subseteq -$\textit{int}$C\}\cap
K(x)=X$\textit{;}

\textit{v) }$K^{-1}:X\rightrightarrows X$\textit{\ is convex valued and for
each }$u\in X,$\textit{\ the set }$\{x\in X:f(x,u)\subseteq -$\textit{int}$%
C\}$\textit{\ is convex}$;$

\textit{Then, there exists }$x^{\ast }\in X$\textit{\ such that }$x^{\ast
}\in K(x^{\ast })$\textit{\ and }$f(x^{\ast },y)\nsubseteq -$\textit{int}$C$%
\textit{\ for each }$y\in K(x^{\ast })$\textit{.\medskip }

\textit{Proof.} Let $P,G:X\rightrightarrows X$ be defined by

$P(x)=\{u\in X:f(x,u)\subseteq -$int$C\},$ for each $x\in X$ and

$G(x)=K(x)\cap P(x),$ for each $x\in X.$

We are going to show that there exists $x^{\ast }\in X$ such that $x^{\ast
}\in K(x^{\ast })$ and $K(x^{\ast })\cap P(x^{\ast })=\emptyset .$

We consider two cases.

Case I.

$A=\{x\in X:$\textit{\ }$P(x)\cap K(x)\neq \emptyset \}=\{x\in X:$\textit{\ }%
$G(x)\neq \emptyset \}$\textit{\ }is nonempty.

The correspondence $G:A\rightrightarrows X,$ defined by $G(x)=K(x)\cap P(x)$
for each $x\in A,$ is nonempty valued on $A.$ We note that for each $u\in X,$
$G^{-1}(u)=P^{-1}(u)\cap K^{-1}(u)$ is a convex set since it is an
intersection of convex sets.

Further, let us define the correspondences $H,L:X\rightrightarrows X$ by

$H(x)=\left\{ 
\begin{array}{c}
G(x),\text{ if }x\in A; \\ 
\emptyset ,\text{\ otherwise}%
\end{array}%
\right. $ and

$L(x)=\left\{ 
\begin{array}{c}
G(x),\text{ \ \ \ \ \ \ \ \ \ \ \ if }x\in A; \\ 
K(x),\text{ \ \ \ \ \ \ \ \ \ otherwise.}%
\end{array}%
\right. $

According to i), $H$ is open valued.

For each $u\in X,$

$H^{-1}(u)=\{x\in X:$\textit{\ }$\mathit{u}\in H(x)\}=$

\ \ \ \ \ \ \ \ \ \ \ $=\{x\in A:$\textit{\ }$\mathit{u}\in G(x)\}=$

\ \ \ \ \ \ \ \ \ \ \ $=\{x\in A:$\textit{\ }$\mathit{u}\in K(x)\cap P(x)\}$

$\ \ \ \ \ \ \ \ =A\cap P^{-1}(u)\cap K^{-1}(u)=$

\ \ \ \ \ \ \ \ \ $=P^{-1}(u)\cap K^{-1}(u).$

$L^{-1}(u)=\{x\in X:$\textit{\ }$\mathit{u}\in L(x)\}=$

\ \ \ \ \ \ \ \ \ \ \ $=\{x\in A:$\textit{\ }$\mathit{u}\in G(x)\}\cup
\{x\in X\setminus A:$\textit{\ }$\mathit{u}\in K(x)\}=$

\ \ \ \ \ \ \ \ \ \ \ $=\{x\in A:$\textit{\ }$\mathit{u}\in K(x)\cap
P(x)\}\cup \{x\in X\setminus A:\mathit{\ u}\in K(x)\}$

$\ \ \ \ \ \ \ \ =(A\cap P^{-1}(u)\cap K^{-1}(u)\cup \lbrack (X\backslash
A)\cap $\textit{\ }$K^{-1}(u)]=$

\ \ \ \ \ \ \ \ \ $=[P^{-1}(u)\cup (X\backslash A)]\cap K^{-1}(u).$

Since for each $u\in X,$ $K^{-1}(u)$ and$P^{-1}(u)$ are convex, then, $%
H^{-1}(u)$ is convex. Therefore, co$H^{-1}(u)\subset L^{-1}(u)$ for each $%
u\in X$.

Assumption iii) implies that there exists $M\in \langle A\rangle $ such that 
$\tbigcup\nolimits_{x\in M}H(x)=X.$

Let us define $F:X\rightrightarrows X$ by $F(x):=X\backslash H(x)$ for each%
\textit{\ }$x\in X.$

Then, $F$ is closed valued and $\tbigcap\nolimits_{x\in M}F(x)=X\backslash
\tbigcup\nolimits_{x\in M}(H(x))=\emptyset .$

According to Lemma 3, we can conclude that $F$ is not a KKM correspondence.
Thus, there exists $N\in \langle X\rangle $ such that\textit{\ }co$%
N\varsubsetneq F(N)=\tbigcup\nolimits_{x\in N}(X\backslash H(x)).$

Hence, there exists $x^{\ast }\in $co$N$ with the property that $x^{\ast
}\in H(x)$ for each $x\in N.$ Therefore, there exists $x^{\ast }\in $co$N$
such that $x^{\ast }\in H(x)$ for each $x\in N,$ which implies $N\subset
H^{-1}(x^{\ast })$. Further, it is true that co$N\subset $co$H^{-1}(x^{\ast
})\subset L^{-1}(x^{\ast })$. Consequently, $x^{\ast }\in $co$N\subset $co$%
H^{-1}(x^{\ast })\subset L^{-1}(x^{\ast }),$ which means that $x^{\ast }\in
L(x^{\ast }),$ that is, $x^{\ast }\in $co$N$ is a fixed point for $L.$

We notice that, according to $ii)$, $x\notin P(x)$ for each $x\in X,$ and
then, $x^{\ast }\notin A.$ This fact is possible since $x^{\ast }\in $co$N$
and assumption iii) asserts that $($co$N\smallsetminus N)\cap A=\emptyset .$
Therefore, $x^{\ast }\in K(x^{\ast })$ and $G(x^{\ast })=K(x^{\ast })\cap
P(x^{\ast })=\emptyset .$

Consequently, there exists $x^{\ast }\in X$ such that $x^{\ast }\in
K(x^{\ast })$ and $f(x^{\ast },u)\nsubseteq -$int$C$ for each $u\in
K(x^{\ast }).$

Case II.

$A=\{x\in X:$\textit{\ }$P(x)\cap K(x)\neq \emptyset \}=\{x\in X:$\textit{\ }%
$G(x)\neq \emptyset \}\mathit{=}\emptyset $.

In this case, $G(x)=\emptyset $ for each $x\in X.$

Let us define $F:X\rightrightarrows X$ by $F(x):=X\backslash K(x)$ for each%
\textit{\ }$x\in X.$ The proof follows the same line as above and we obtain
that there exists $x^{\ast }\in X$ such that $x^{\ast }\in K(x^{\ast }).$

Obviously, $G(x^{\ast })=\emptyset .$ Consequently, the conclusion holds in
Case II.

\begin{remark}
Assumption i) can be replaced with
\end{remark}

i') $K$\ is closed valued and for each $x\in X,$ the set $\{u\in
X:f(x,u)\subseteq -$int$C\}$ is closed.\textit{\ }

In this case, $F$ is open valued.\smallskip

Theorem 7 can be stated as follows.

\begin{theorem}
\textit{Let }$Z$\textit{\ be a Hausdorff topological vector space and let }$%
X $\textit{\ be a nonempty compact convex subset of a\ topological vector
space }$E$\textit{. Let }$C\subset Z$\textit{\ be a pointed closed convex
cone with nonempty interior and let }$K:X\rightrightarrows X$\textit{\ and }$%
f:X\times X\rightrightarrows Z$\textit{\ be correspondences with nonempty
values. Assume that:}
\end{theorem}

\textit{i) }$K$\textit{\ is open (resp. closed) valued and for each }$x\in
X, $\textit{\ the set }$\{u\in X:f(x,u)\cap $\textit{int}$C=\emptyset \}$%
\textit{\ is open (resp.closed); }

\textit{ii) \ }$f(x,x)\cap $\textit{int}$C=\emptyset $\textit{\ for each }$%
x\in X$\textit{; }

\textit{iii) if }$A=\{x\in X:$\textit{\ there exists }$u\in K(x)$\textit{\
such that }$f(x,u)\cap $\textit{int}$C=\emptyset \},$\textit{\ then for each 
}$N\in \langle A\rangle ,$\textit{\ }$($\textit{co}$N\smallsetminus N)\cap
A=\emptyset ;$

\textit{iv) there exists }$M\in \langle A\rangle $\textit{\ such that }$%
\bigcup\nolimits_{x\in M}\{u\in X:f(x,u)\cap $\textit{int}$C=\emptyset
\}\cap K(x)=X$\textit{;}

\textit{v) }$K^{-1}:X\rightrightarrows X$\textit{\ is convex valued and for
each }$u\in X,$\textit{\ the set }$\{x\in X:f(x,u)\cap $\textit{int}$%
C=\emptyset \}$\textit{\ is convex}$.$

\textit{Then, there exists }$x^{\ast }\in X$\textit{\ such that }$x^{\ast
}\in K(x^{\ast })$\textit{\ and }$f(x^{\ast },y)\cap $\textit{int}$C\neq
\emptyset $\textit{\ for each }$y\in K(x^{\ast })$\textit{.\medskip }

\textit{Proof.} Let $P,G:X\rightrightarrows X$ be defined by

$P(x)=\{u\in X:f(x,u)\cap $int$C=\emptyset \},$ for each $x\in X$ and

$G(x)=K(x)\cap P(x),$ for each $x\in X.$

We are going to show that there exists $x^{\ast }\in X$ such that $x^{\ast
}\in K(x^{\ast })$ and $K(x^{\ast })\cap P(x^{\ast })=\emptyset .$

The rest of the proof follows a similar line as the proof of Theorem
6.\smallskip

Now, we are proving the existence of solutions for a general vector
equilibrium problem concerning correspondences under new assumptions.

\begin{theorem}
\textit{Let }$Z$\textit{\ be a Hausdorff topological vector space and let }$%
X $\textit{\ be a nonempty compact convex subset of a\ topological vector
space }$E$\textit{. Let }$C\subset Z$\textit{\ be a pointed closed convex
cone with nonempty interior and let }$K:X\rightrightarrows X$\textit{\ and }$%
f:X\times X\rightrightarrows Z$\textit{\ be correspondences with nonempty
values. Assume that:}
\end{theorem}

\textit{i) }$K$\textit{\ is open valued and for each }$(x,y)\in X\times X$%
\textit{\ with the property that }$f(x,y)\subseteq -$\textit{int}$C,$\textit{%
\ there exists }$z=z_{x,y}\in X$\textit{\ such that }$y\in $\textit{int}$%
_{X}\{u\in X:f(z_{x,y},u)\subseteq -$int$C\}$\textit{; }

\textit{ii) }$f(x,x)\nsubseteq -$\textit{int}$C$\textit{\ for each }$x\in X$%
\textit{; }

\textit{iii) there exists }$M\in \langle X\rangle $\textit{\ such that }$%
\bigcup\nolimits_{x\in M\cap A}[\bigcup\nolimits_{y\in K(x),f(x,y)\subseteq -%
\text{int}C}($\textit{int}$_{X}\{u\in X:f(z_{x,y},u)\subseteq -$\textit{int}$%
C\}\cap K(x)]\bigcup \bigcup\nolimits_{x\in M\backslash A}K(x)=X,$\textit{\
where }$A=\{x\in X:$\textit{\ there exists }$u\in K(x)$\textit{\ such that }$%
f(x,u)\subseteq -$\textit{int}$C\}$\textit{;}

\textit{iv) for each }$u\in X,$\textit{\ the set }$\{x\in X:f(x,u)\subseteq
- $\textit{int}$C\}\cup \lbrack (X\backslash A)\cap (K^{-1}(u)]$\textit{\ is
convex}$;$

\textit{v) for each }$x\in A,$\textit{\ }$\bigcup\nolimits_{y\in
K(x),f(x,y)\subseteq -\text{int}C}$\textit{int}$_{X}\{u\in
X:f(z_{x,y},u)\subseteq -$\textit{int}$C\}\subset \{u\in X:f(x,u)\subseteq -$%
int$C\}$\textit{;}

\textit{Then, there exists }$x^{\ast }\in X$\textit{\ such that }$x^{\ast
}\in K(x^{\ast })$\textit{\ and }$f(x^{\ast },y)\nsubseteq -$\textit{int}$C$%
\textit{\ for each }$y\in K(x^{\ast })$\textit{.\medskip }

\textit{Proof.} Let $P,G:X\rightrightarrows X$ be defined by

$P(x)=\{u\in X:f(x,u)\subseteq -$int$C\},$ for each $x\in X$ and

$G(x)=K(x)\cap P(x),$ for each $x\in X.$

We are going to show that there exists $x^{\ast }\in X$ such that $x^{\ast
}\in K(x^{\ast })$ and $K(x^{\ast })\cap P(x^{\ast })=\emptyset .$

We consider two cases.

Case I.

$A=\{x\in X:$\textit{\ }$P(x)\cap K(x)\neq \emptyset \}=\{x\in X:$\textit{\ }%
$G(x)\neq \emptyset \}$\textit{\ }is nonempty.

The correspondence $G:A\rightrightarrows X,$ defined by $G(x)=K(x)\cap P(x)$
for each $x\in A,$ is nonempty valued on $A$.

Further, let us define the correspondences $H,L:X\rightrightarrows X$ by

$H(x)=\left\{ 
\begin{array}{c}
G(x),\text{ if }x\in A; \\ 
K(x),\text{\ otherwise}%
\end{array}%
\right. $ and

$L(x)=\left\{ 
\begin{array}{c}
P(x),\text{ if }x\in A; \\ 
K(x),\text{\ otherwise}%
\end{array}%
\right. $

For each $u\in X,$

$L^{-1}(u)=\{x\in X:$\textit{\ }$\mathit{u}\in L(x)\}=$

\ \ \ \ \ \ \ \ \ \ \ $=\{x\in A:$\textit{\ }$\mathit{u}\in P(x)\}\cup
\{x\in X\setminus A:$\textit{\ }$\mathit{u}\in K(x)\}=$

$\ \ \ \ \ \ \ \ =(A\cap P^{-1}(u))\cup \lbrack (X\backslash A)\cap $\textit{%
\ }$K^{-1}(u)]=$

\ \ \ \ \ \ \ \ \ $=P^{-1}(u)\cup \lbrack (X\backslash A)\cap K^{-1}(u)].$

According to i) $G$ is transfer-open valued. Assumptions i) and iii) imply
that there exists $M\in \langle X\rangle $ and for each $x\in M$ and $y\in
H(x),$ there exists $z_{x,y}\in X$ such that $y\in $int$_{X}H(z_{x,y})\cap
H(x)$ and $\tbigcup\nolimits_{x\in M}(\tbigcup\nolimits_{y\in H(x)}$int$%
_{X}H(z_{x,y}))=X$. In addition, assumption v) implies $\tbigcup\nolimits_{y%
\in H(x)}$int$_{X}H(z_{x,y})\subseteq P(x)$ for each $x\in A.$ We note that
if $x\in X\backslash A,$ then $H(x)=K(x)$ is open and $y\in H(x)$ implies $%
z_{x,y\text{ }}=x$ and $y\in H(z_{x.y})=$int$H(x).$ In this last case it is
obvious that $\tbigcup\nolimits_{y\in H(x)}$int$_{X}H(z_{x,y})=\tbigcup%
\nolimits_{y\in H(x)}$int$_{X}H(x)=\tbigcup\nolimits_{y\in H(x)}H(x)=H(x).$

Let us define $F:X\rightrightarrows X$ by $F(x):=X\backslash
\tbigcup\nolimits_{y\in H(x)}$int$_{X}H(z_{x,y})$ for each\textit{\ }$x\in
H. $

Then, $F$ is closed valued and $\tbigcap\nolimits_{x\in M}F(x)=X\backslash
\tbigcup\nolimits_{x\in M}(\tbigcup\nolimits_{y\in H(x)}$int$%
_{X}H(z_{x,y}))=\emptyset .$

According to Lemma 3, we can conclude that $F$ is not a KKM correspondence.
Thus, there exists $N\in \langle X\rangle $ such that\textit{\ }co$%
N\varsubsetneq F(N)=\tbigcup\nolimits_{x\in N}(X\backslash
\tbigcup\nolimits_{y\in H(x)}$int$_{X}H(z_{x,y})).$

Hence, there exists $x^{\ast }\in $co$N$ with the property that $x^{\ast
}\in \tbigcup\nolimits_{y\in H(x)}$int$_{X}H(z_{x,y})$ for each $x\in N.$ If 
$x\in A,$ $\tbigcup\nolimits_{y\in H(x)}$int$_{X}H(z_{x,y})\subset P(x)$ and
if $x\in X\backslash A,$ $\tbigcup\nolimits_{y\in H(x)}$int$%
_{X}H(z_{x,y})=H(x).$ Therefore, there exists $x^{\ast }\in $co$N$ such that 
$x^{\ast }\in L(x)$ for each $x\in N,$ which implies $N\subset
L^{-1}(x^{\ast })$. Further, it is true that co$N\subset $co$L^{-1}(x^{\ast
})=L^{-1}(x^{\ast })$. Consequently, $x^{\ast }\in $co$N\subset $co$%
L^{-1}(x^{\ast })=L^{-1}(x^{\ast }),$ which means that $x^{\ast }\in
L(x^{\ast }).$ We notice that, according to $ii)$, $x\notin P(x)$ for each $%
x\in X,$ and then, $x^{\ast }\notin A.$ Therefore, $x^{\ast }\in K(x^{\ast
}) $ and $G(x^{\ast })=K(x^{\ast })\cap P(x^{\ast })=\emptyset .$

Consequently, there exists $x^{\ast }\in X$ such that $x^{\ast }\in
K(x^{\ast })$ and $f(x^{\ast },u)\nsubseteq -$int$C$ for each $u\in
K(x^{\ast }).$

Case II.

$A=\{x\in X:$\textit{\ }$P(x)\cap K(x)\neq \emptyset \}=\{x\in X:$\textit{\ }%
$G(x)\neq \emptyset \}\mathit{=}\emptyset $.

In this case, $G(x)=\emptyset $ for each $x\in X.$ Let us define $%
F:X\rightrightarrows X$ by $F(x):=X\backslash K(x)$ for each\textit{\ }$x\in
X.$ The proof follows the same line as above and we obtain that there exists 
$x^{\ast }\in X$ such that $x^{\ast }\in K(x^{\ast }).$

Obviously, $G(x^{\ast })=\emptyset .$ Consequently, the conclusion holds in
Case II.\smallskip

Now, we are establishing Theorem 9.

\begin{theorem}
\textit{Let }$Z$\textit{\ be a Hausdorff topological vector space and let }$%
X $\textit{\ be a nonempty compact convex subset of a\ topological vector
space }$E$\textit{. Let }$C\subset Z$\textit{\ be a pointed closed convex
cone with nonempty interior and let }$K:X\rightrightarrows X$\textit{\ and }$%
f:X\times X\rightrightarrows Z$\textit{\ be correspondences with nonempty
values. Assume that:}
\end{theorem}

\textit{i) }$K$\textit{\ is open valued and for each }$(x,y)\in X\times X$%
\textit{\ with the property that }$f(x,y)\cap $\textit{int}$C=\emptyset ,$%
\textit{\ there exists }$z=z_{x,y}\in X$\textit{\ such that }$y\in $\textit{%
int}$_{X}\{u\in X:f(z_{x,y},u)\cap $\textit{int}$C=\emptyset \}$\textit{; }

\textit{ii) }$f(x,x)\cap $\textit{int}$C=\emptyset $\textit{\ for each }$%
x\in X$\textit{; }

\textit{iii) there exists }$M\in \langle X\rangle $\textit{\ such that }$%
\bigcup\nolimits_{x\in M\cap A}[\bigcup\nolimits_{y\in K(x),f(x,y)\cap \text{%
int}C=\emptyset }($\textit{int}$_{X}\{u\in X:f(z_{x,y},u)\cap $\textit{int}$%
C=\emptyset \}\cap K(x)]\bigcup \bigcup\nolimits_{x\in M\backslash A}K(x)=X,$%
\textit{\ where }$A=\{x\in X:$\textit{\ there exists }$u\in K(x)$\textit{\
such that }$f(x,u)\cap $\textit{int}$C=\emptyset \}$\textit{;}

\textit{iv) for each }$u\in X,$\textit{\ the set }$\{x\in X:f(x,u)\cap $%
\textit{int}$C=\emptyset \}\cup \lbrack (X\backslash A)\cap (K^{-1}(u)]$%
\textit{\ is convex}$;$

\textit{v) for each }$x\in A,$\textit{\ }$\bigcup\nolimits_{y\in
K(x),f(x,y)\cap \text{int}C=\emptyset }$\textit{int}$_{X}\{u\in
X:f(z_{x,y},u)\cap $\textit{int}$C=\emptyset \}\subset \{u\in
X:f(x,u)\subseteq -$int$C\}$\textit{;}

\textit{Then, there exists }$x^{\ast }\in X$\textit{\ such that }$x^{\ast
}\in K(x^{\ast })$\textit{\ and }$f(x^{\ast },y)\cap $\textit{int}$C\neq
\emptyset $\textit{\ for each }$y\in K(x^{\ast })$\textit{.\medskip }

\section{CONCLUDING\ REMARKS}

In the first part of this paper, we have proved the existence of solutions
for generalized quasi-equilibrium problems on Banach spaces. The involved
correspondences are almost lower semicontinuous. Then, we have established a
new result in topological vector spaces, in the case of lower semicontinuity
of the correspondences. Our research extends on some results which exist in
literature and is based on fixed point theorems. In the second part, we have
obtained new equilibrium theorems by applying the KKM principle. This study
will be continued by considering abstract convex spaces and generalized KKM
theorems.

\end{document}